\documentclass[twocolumn]{autart}    

\usepackage{graphicx}          
\usepackage{amsmath}
\usepackage{amsfonts}
\usepackage{xcolor}
\usepackage{url}

\usepackage{multicol}

\usepackage{algorithm}
\usepackage{algpseudocode}

\usepackage{cite}
\usepackage{enumitem}

\usepackage{graphicx}
\usepackage{caption}
\usepackage{subcaption}

\usepackage{multirow}
\usepackage{hyperref}

\begin{document}

\begin{frontmatter}

\title{Semismooth Newton Methods for\\ Risk-Averse Markov Decision Processes\thanksref{footnoteinfo}} 

\thanks[footnoteinfo]{This paper was not presented at any IFAC 
meeting. Corresponding author Matilde Gargiani.}

\author[Zurich]{Matilde Gargiani}\ead{gmatilde@ethz.ch},    
\author[Zurich]{Francesco Micheli}\ead{frmicheli@ethz.ch},               
\author[Zurich]{Anastasios Tsiamis}\ead{atsiamis@ethz.ch},
\author[Zurich]{John Lygeros}\ead{jlygeros@ethz.ch}  

\address[Zurich]{Automatic Control Laboratory, Physikstrasse 3, 8092 Zurich.}  
          

\begin{abstract}                          
Inspired by semismooth Newton methods, we propose a general framework for designing solution methods with convergence guarantees for risk-averse Markov decision processes. Our approach accommodates a wide variety of risk measures by leveraging the assumption of Markovian coherent risk measures. To demonstrate the versatility and effectiveness of this framework, we design three distinct solution methods, each with proven convergence guarantees and competitive empirical performance.
Validation results on benchmark problems demonstrate the competitive performance of our methods. 
Furthermore, we establish that risk-averse policy iteration can be interpreted as an instance of semismooth Newton's method. This insight explains its superior convergence properties compared to risk-averse value iteration. The core contribution of our work, however, lies in developing an algorithmic framework inspired by semismooth Newton methods, rather than evaluating specific risk measures or advocating for risk-averse approaches over risk-neutral ones in particular applications.
\end{abstract}

\end{frontmatter}

\section{Introduction}
Markov Decision Processes (MDPs) are a very powerful mathematical tool that can be used to model real-world problems arising in a variety of different fields~\cite{mdp_applications}, from finance~\cite{mdp_finance} to epidemiology~\cite{mdp_epidemiology}. Optimal decision making within an MDP framework is generally tackled with Dynamic Programming (DP), initiated in the 50s by R. Bellman~\cite{Bellman1957}. Classical DP only deals with expected performance criteria. These often are not appropriate metrics to measure performance as they do not take into account the potential fluctuations of the cost around the mean nor its sensitivity to modeling errors. In many applications variability from the expected value detrimentally impacts the performance of the system, such as in finance~\cite{coherent_riskmeasures} and some robotics applications~\cite{pavone_robotics}. For instance, the authors in~\cite{pavone_robotics} bring up as example an autonomous navigation system, where a robot attempts to minimize the expected length of the traveled path. This objective will likely push the robot to travel very close to obstacles, so a deviation from the planned path, due, for instance, to modeling errors, may result in a collision that damages the robot. Likewise,~\cite[Example 1]{risk_averse_DP} provides a simple example where stochastic realizations with limited deviations from the expected value are preferable even if they entail a higher cost. Risk-averse DP replaces the expected value with more general risk-measures, which provide a way to both quantify the risk and to ensure robustness to modeling errors. This variation of the classical DP setting comes with some technical challenges; for instance, the problem may lose its Markovian structure breaking the recursive form of the Bellman equations. This creates theoretical and computational hurdles~\cite{MANNOR2013645} and requires solution methods that are designed ad-hoc for specific risk-measures~\cite{RiskSensitiveAR, PG_CvAR}. These methods are generally based on approximations of the original problem that can be efficiently solved, but in general do not retrieve the optimal policy~\cite{pavone_robotics}.
Owing to time inconsistency, it can indeed be difficult or impossible to construct DP formulations and to derive the corresponding Bellman operators~\cite{MCR_PG_proof, risk_averse_DP}. 

In~\cite{risk_averse_DP} Ruszczy\'nski introduces the concept of Markovian coherent risk measure and derives a risk-averse extension of the classical Bellman equations for discounted infinite-horizon MDPs. This leads to the definition of risk-averse variants of value iteration and policy iteration, that share similar convergence properties as those of their risk-neutral counterparts.
Many researchers have adopted Ruszczy\'nski's framework and use Markov coherent risk measures~\cite{risk_averse_DP} as a time-consistent surrogate which satisfies Bellman's principle of optimality for risk-averse MDPs~\cite{MCR_PG_proof, Tamar}. This makes it possible to tackle risk-averse MDPs with a unified framework and deploy the same techniques and logic of risk-neutral dynamic programming.

In this work, we propose to adopt a new perspective to tackle the solution of the risk-averse extension of the Bellman equations proposed by Ruszczy\'nski. In particular, in the risk-neutral scenario, the connection between policy iteration and semismooth Newton methods has lead to significant theoretical and algorithmic insights~\cite{puterman, gargiani_newton, qi1997semismooth} and also the design of novel efficient methods to tackle large-scale MDPs~\cite{gargiani_newton, gargiani_igmres, gargiani_ipi,qi1997semismooth}. Inspired by these results, we study the deployment of semismooth Newton methods for risk-averse MDPs with Markovian coherent risk-measures. For that, we start from the framework proposed by Ruszczy\'nski and focus on the class of discounted infinite-horizon MDPs with finite state and action spaces. For this class of problems
\begin{enumerate}
    \item we propose a new framework based on semismooth Newton methods to design novel solution methods for risk-averse MDPs with Markovian risk-measures,
    \item we deploy our framework to design three different methods to solve risk-averse DP problems and provide a mathematical characterization of their convergence properties,
    \item we prove that the risk-averse policy iteration method designed by Ruszczy\'nski in~\cite{risk_averse_DP} is itself an instance of a semismooth Newton method. This connection helps to shed some light on the superior convergence properties of this method with respect to risk-averse value iteration,
    \item taking inspiration from risk-neutral DP methods, we propose an inexact variant of our methods and prove its global convergence. This variant is the risk-averse analog of optimistic policy iteration.
    \item we define sufficient local structural properties which guarantee fast quadratic local convergence for methods developed within this framework.
\end{enumerate}

We exemplify the practical utility of our framework by considering Conditional Value at Risk (CVaR) MDPs~\cite{RiskSensitiveAR} and we study the empirical performance of the proposed methods, including our risk-averse version of optimistic policy iteration, on this class of risk-averse MDPs.
It is important to note that the paper’s contribution does not lie in promoting the importance of adopting a risk-averse approach or in comparing various risk measures. Instead, it focuses on the design and analysis of solution methods to effectively address this class of problems in a unified way.

The paper is organized as follows: in Section~\ref{sec:problem_setting} we introduce the problem setting and review the existing methods for risk-averse DP. In Section~\ref{sec:newton_f} we discuss Newton approximations schemes and illustrate how to deploy them within a semismooth Newton method. We also introduce the risk-averse Bellman residual function and then propose three different Newton approximation schemes for it. In Sections~\ref{sec:SNM1}-\ref{sec:SNM3} we use the proposed Newton approximation schemes to design three different instances of the semismooth Newton method for the risk-averse Bellman residual function. For each method we mathematically study its convergence properties and consider the computational aspects. Finally, in Section~\ref{sec:numerical-experiments}, we study the empirical performance of the proposed methods on CVaR MDPs, demonstrating their practical utility in solving risk-sensitive applications.
\section*{Notation}
We use $\mathbb{S}^n$ and $\mathbb{S}^{n\times n}$ to denote the simplex of dimension $n$, \textit{i.e.} $\left\{s\in\mathbb{R}^n\,|\, \sum_{i=1}^n s_i = 1,\,s_i\geq 0 \,\,\,\forall\, i\in\left\{1,\dots,n \right\} \right\}$, and set of row-stochastic matrices of dimension $n\times n$, respectively.
Given an operator $\mathcal{O}:\mathbb{R}^n \rightarrow \mathbb{R}^n$, an integer $m>0$ and $v\in \mathbb{R}^n$, we denote with $\left(\mathcal{O}\right)^m v$ the application of the $\mathcal{O}$-operator $m$-times to $v$. We denote with $v(i)$ the $i$-th entry of $v$. Given a matrix $A \in \mathbb{R}^{n\times m}$ we use $[A]_{i, :}$ to denote the $i$-th row, with $i\in \left\{1, \dots, n \right\}$. We use $e$ to denote the vector of ones in $\mathbb{R}^n$ and $(t)_{+} = \max\left\{0,\,t \right\}$ where $t\in\mathbb{R}^n$ and the max is understood element-wise. We use $\odot$ to denote the element-wise product. We use $\mathcal{U}_x$ to denote a neighborhood of $x\in \mathbb{R}^n$.
For $s,\,\tilde{s}\in\mathbb{R}^n$ we intend $s\leq \tilde{s}$ element-wise. Let $a,\,b\in\mathbb{R}$ with $a<b$, we use $[a,b]^n$ to denote the set $\underbrace{[a,b]\times\dots\times [a,b]}_{n}$.

\section{Risk-Averse Markov Decision Processes}\label{sec:problem_setting}
We consider discounted stationary MDPs over an infinite-horizon with finite state and action spaces. The MDPs in this class can be mathematically described by 5-elements $\left\{\mathcal{S},\,\mathcal{A},\,\gamma\,, P,\, g \right\}$, where, without loss of generality, $\mathcal{S} = \left\{1,\dots, n \right\}$ and $\mathcal{A} = \left\{1,\dots, m \right\}$ are the state and action spaces, respectively; $\gamma\in(0,1)$ is the discount factor; $P:\mathcal{S}\times \mathcal{A}  \rightarrow \mathbb{S}^n$ is the controlled kernel and $g:\mathcal{S}\times \mathcal{A} \rightarrow [-R, \,R]$ with $R>0$ is the stage-cost function. In addition, we denote with $\mathcal{A}(x)\subseteq \mathcal{A}$ the subset of actions which are allowed in state $x\in\mathcal{X}$. An admissible stationary and deterministic policy, simply called policy from now on, is a function  $\pi:\mathcal{X} \rightarrow \mathcal{U}$ that maps states into actions, with $\pi(x)\in \mathcal{A}(x)$ for all $x\in\mathcal{X}$. We define with $\Pi$ the set of all policies.
The aim of classical dynamic programming in the considered problem setting is to find a policy $\pi\in\Pi$ that minimizes the following objective function for any $i\in \mathcal{S}$
\begin{equation}\label{eq:risk-neutral_DP_of}
    v^{\pi}(i) =\lim_{H\rightarrow \infty} \mathbb{E}\left[ \sum_{t=0}^{H-1} \gamma^t g(s_t, \pi(s_t)) \,\Big\vert\, s_0 = i\right]\,.
\end{equation}
Classical DP is therefore concerned with expected performance metrics~\cite{bertsekas_book}. Since in many problems of practical interest expected criteria are not appropriate to measure performance, in this work we look into replacing the expected value operator in~\eqref{eq:risk-neutral_DP_of} by more general risk-measures to model risk aversion. In particular, we aim at finding the policy $\pi\in \Pi$ which minimizes the following objective function for any $i\in \mathcal{S}$
\begin{equation}\label{eq:obj_fun}
    v_{\mu}^{\pi}(s_0) = \lim_{H\rightarrow \infty} \mu \left( \sum_{t = 0}^{H-1} \gamma^t g(s_t, \pi(s_t)) \right)\,,
\end{equation}
with $s_0 = i$ and $\mu:\mathcal{L}\rightarrow \mathbb{R}$ a risk-measure, where $\mathcal{L}$ is the space of all random variables from the canonical sample space $\Omega$ to the space of trajectories. From now on, with a slight abuse of notation we use $v^{\pi}$ in place of $v_{\mu}^{\pi}$. Equation~\eqref{eq:obj_fun} defines the risk-averse cost associated to policy $\pi$.
We call optimal cost and - with a slight abuse of notation - denote it with $v^*$ the minimum of~\eqref{eq:obj_fun} with respect to all $\pi\in \Pi$. Notice that the minimum is always attained in this setting since the set of policies is finite. Finally, we define as optimal and denote with $\pi^*$ a policy $\pi\in\Pi$ that attains the optimal cost. 

We consider the case where $\mu$ is a stationary Markov risk-measure in the sense of Definition 6 in~\cite{risk_averse_DP}.
Consequently, considering the MDP defined by the tuple $\left\{\mathcal{S},\,\mathcal{A},\,\gamma,\, P,\,g \right\}$, there exists a risk-transition mapping $\chi:\mathbb{R}^n \times \mathbb{S}^n\rightarrow \mathbb{R}$~\cite[Definition 5]{risk_averse_DP} such that for all $v\in \mathbb{R}^n$ and $a_t \in \mathcal{A}(x_t)$
\begin{equation}\label{eq:sigma_equivalence}
    \mu(v(s_{t+1})) = \chi(v, P(s_t,\,a_t))\,.
\end{equation}
The Markovian assumption on $\mu$ and relation~\eqref{eq:sigma_equivalence} are fundamental for our formulation as otherwise we can not expect to obtain a Markovian optimal policy for~\eqref{eq:obj_fun} since general risk-measures could depend on the entire history of the process~\cite{risk_averse_DP}. The Markovian assumption allows us to re-write~\eqref{eq:obj_fun} recursively as 
\begin{equation}\label{eq:markovian_rho}
\begin{aligned}
    v^{\pi}(s_0) = \lim_{H\rightarrow \infty}  g(s_0, &\pi(s_0))  + \gamma \mu(g(s_1, \pi(s_1))\\& + \gamma \mu(g(s_2, \pi(s_2)) + \dots ))\,.
\end{aligned}
\end{equation}
We refer to~\cite{risk_averse_DP} for more details and a formal proof on how to unroll~\eqref{eq:obj_fun} into~\eqref{eq:markovian_rho}.

Following~\cite{risk_averse_DP}, given the controlled kernel $P$ and for any $(i,a)\in\mathcal{S}\times \mathcal{A}$, in our setting there always exists a closed, bounded and compact convex set of distributions $\mathcal{P}_{\chi}(P(i,a)) \subset \mathbb{S}^n$ such that for all $v\in \mathbb{R}^n$ we have that
\begin{equation}\label{eq:dual_rewrite}
    \chi(v, P(i,a)) = \max_{\tilde{p} \in \mathcal{P}_{\chi}(P(i, a))} \langle v, \tilde{p} \rangle\,.
\end{equation}
This reformulation of the risk-transition map has an interesting distributionally robust interpretation: the right-hand side in~\eqref{eq:dual_rewrite} can be rewritten as the expected value of the random variable $v(s_{t+1})$ with respect to the worst-case distribution in a set of perturbed distributions that is determined by the specific choice of $\chi$ and the controlled kernel of the underlying MDP. We can therefore use this setting when the controlled kernel is subject to uncertainty and use the risk-transition map to model this uncertainty in a distributionally robust fashion.   


As in the risk-neutral setting, it is convenient to define two operators which we call \textit{risk-averse Bellman operators}. Given a policy $\pi:\mathcal{S}\rightarrow \mathcal{A}$, we define $D_{\pi}:\mathbb{R}^n \rightarrow \mathbb{R}^n$ as
\begin{equation}\label{eq:risk-averse_BE_pi}
    (D_{\pi} v)(i) = g(i, \pi(i)) + \gamma \chi(v, P(i, \pi(i))) \quad i\in\mathcal{S}\,,
\end{equation}
and we call it the risk-averse Bellman operator associated to $\pi$.
We also define the risk-averse Bellman operator for optimality $D : \mathbb{R}^n \rightarrow \mathbb{R}^n$ as
\begin{equation}\label{eq:risk-averse_BE}
    (D v)(i) = \min_{a \in \mathcal{A}(i)}\left\{ g(i, a) + \gamma \chi(v, P(i, a)) \right\} \quad i\in \mathcal{S}\,.
\end{equation}
Notice that, unlike the risk-neutral scenario,~\eqref{eq:risk-averse_BE_pi} is not necessarily linear nor smooth because of the potential non-linearity and non-smoothness of $\chi$. 
It can be shown that the risk-averse Bellman operators~\eqref{eq:risk-averse_BE_pi} and~\eqref{eq:risk-averse_BE} are monotone~\cite[Lemma 1]{risk_averse_DP} and $\gamma$-contractive in the infinity norm~\cite[Lemma 2]{risk_averse_DP}. In addition,~\eqref{eq:obj_fun} has a unique minimum $v^*$ which is also the unique fixed-point of the $D$-operator~\cite[Theorem 4]{risk_averse_DP}. Finally, given a policy $\pi\in\Pi$,~\eqref{eq:obj_fun} is the unique fixed-point of the $D_{\pi}$-operator~\cite[Lemma 4]{risk_averse_DP}.
We summarize these properties in the following lemma.
\begin{lem}\label{lemma:fundamental_properties_D}
    The risk-averse Bellman operators $D$ and $D_{\pi}$ are $\gamma$-contractive in the infinite-norm and monotone. In addition, they have $v^*$ and $v^{\pi}$ as their unique fixed-points, respectively.
\end{lem} 
\noindent\textit{Proof.} The proof follows from~\cite{risk_averse_DP}, since the required continuity assumptions are always verified for finite MDPs.

In the following lemma we show that, as their risk-neutral counterparts~\cite[Lemma 1.1.2]{bertsekas_book}, also the risk-averse Bellman operators enjoy the constant shift property. This will be used in Section~\ref{sec:SNM2} to prove that the sequence of iterates generated by Algorithm~\ref{alg:riskaverse_opi} is globally convergent to $v^*$.

\begin{lem}\label{lm:constant_shift_D}
 For any $v\in \mathbb{R}^n$, $b\in\mathbb{R}$, stationary policy $\pi\in\Pi$ and for $k=0,1,\dots,$ then
        \begin{equation}\label{eq:shift_delta}
            \left(D\right)^{k}(v + b e) = \left(D\right)^{k}v + \gamma^{k} b e\,, 
        \end{equation}
        \begin{equation}\label{eq:shift_deltapi}
            \left(D_{\pi}\right)^{k}(v + b e) = \left(D_{\pi}\right)^{k} v + \gamma^{k} b e\,.
        \end{equation}
\end{lem}
\noindent \textit{Proof.} To simplify the notation, in the following we use $\mathcal{P}_{\chi}$ in place of $\mathcal{P}_{\chi}(P(i,\,\pi(i)))$. For all $i\in\mathcal{S}$
\begin{equation}
    \begin{aligned}
        &(D(v + be))(i) = \min_{\pi\in\Pi} \left\{ g(i,{\pi}(i)) + \gamma \max_{\tilde{p} \in \mathcal{P}_{\chi}}\, \langle \tilde{p}, v + be \rangle  \right\}\\
        &= \min_{\pi\in\Pi} \left\{ g(i,{\pi}(i)) + \gamma \max_{\tilde{p} \in \mathcal{P}_{\chi}} \sum_{j\in\mathcal{S}} \tilde{p}(j) \cdot (v(j) + b)  \right\}\\
        &\overset{(a)}{=} \min_{\pi\in\Pi} \left\{ g(i,{\pi}(i)) + \gamma \max_{\tilde{p} \in \mathcal{P}_{\chi}}\, \langle \tilde{p} , v \rangle + \gamma b \right\}\\
        &= (D v )(i) + \gamma b\,,
    \end{aligned}
\end{equation}
where $(a)$ follows from the fact that $\mathcal{P}_{\chi}\subset \mathbb{S}^n$.
Repeating this argument $k$-times leads to~\eqref{eq:shift_delta}.
Similar step prove~\eqref{eq:shift_deltapi}.~\qed

Making use of~\eqref{eq:risk-averse_BE} and the results of Theorem 4 in~\cite{risk_averse_DP}, we can equivalently and compactly re-write the minimization of~\eqref{eq:obj_fun} over the set of policies as follows
\begin{equation}\label{eq:compact_riskaverse_BE}
    v = Dv\,,
\end{equation}
which is verified if and only if $v = v^*$.
As in the risk-neutral scenario, dynamic programming comprises solution methods for~\eqref{eq:compact_riskaverse_BE}. From the Banach theorem~\cite{Rockafellar2000OptimizationOC} it directly follows that the corresponding Picard-Banach iteration of $v$ linearly converges to $v^*$ from any $v_0\in\mathbb{R}^n$. This method, whose $(k+1)$-iteration for $k\geq 0$ is
\begin{equation}
    v_{k+1} = {D}v_k\,,
\end{equation}
is also known as risk-averse value iteration (risk-averse VI)~\cite{risk_averse_DP}. Likewise, given a policy $\pi\in\Pi$, we call the Picard-Banach iteration of $D_{\pi}$ risk-averse VI for policy evaluation. As for risks averse VI, global linear convergence in the infinity norm follows from the $\gamma$ contraction property of Lemma~\ref{lemma:fundamental_properties_D}.

An alternative method to solve~\eqref{eq:compact_riskaverse_BE} is the so-called risk-averse policy iteration method, which, as suggested by the name, generalizes the risk-neutral policy iteration algorithm~\cite[Chapter 2]{bertsekas_book} to the risk-averse setting.  This scheme is proposed in~\cite{risk_averse_DP}, by replacing the risk-neutral Bellman operators with their risk-averse extensions. In particular, starting from $v_0 \in \mathbb{R}^n$, for $k=0,1,2,\dots$ the algorithm alternates a policy improvement step where $\pi_{k+1}$ is selected according to
\begin{equation}
    D_{\pi_{k+1}} v_k = D v_k\,,
\end{equation}
and a policy evaluation step where the new iterate $v_{k+1}$ is computed as
\begin{equation}\label{eq:policy_evaluation}
    v_{k+1} = {D}_{\pi_{k+1}} v_{k+1}\,.
\end{equation}
Notice that Equation~\eqref{eq:policy_evaluation} requires the computation of the fixed-point of the ${D}_{\pi_{k+1}}$-operator.
The $(k+1)$-iterate is therefore the risk-averse cost associated to the current greedy-policy $\pi_{k+1}$. Unlike risk-neutral policy iteration, the policy evaluation step in~\eqref{eq:policy_evaluation} requires the solution of a potentially non-smooth and non-linear system of equations. Assuming that at every iteration~\eqref{eq:policy_evaluation} is solved exactly, Ruszczyński in~\cite[Theorem 6]{risk_averse_DP} proves that the generated sequence of iterates is non-increasing and globally converges to $v^*$ with a rate $\gamma$ in the infinity-norm.

\section{Semismooth Newton Methods}\label{sec:newton_f}
We start by re-writing~\eqref{eq:risk-averse_BE} as a root finding problem
\begin{equation}\label{eq:risk-averse_rootfinding}
    \underbrace{v - Dv}_{r(v)} = 0\,,
\end{equation}
where $r:\mathbb{R}^n\rightarrow \mathbb{R}^n$ is non-smooth and non-linear because of the non-linearity and non-smoothness of the $D$-operator. In analogy with the risk-neutral scenario, we call $r$ the \textit{risk-averse Bellman residual function}. Because of the underlying assumptions on the risk-measure in~\eqref{eq:obj_fun}, the risk-averse Bellman residual function has $v^*$ as unique root.  

The workhorse methods in mathematical analysis to tackle general non-smooth root finding problems are semismooth Newton methods and their variants~\cite{facchinei_book, izmailov_book}. These are iterative methods that replace the residual function with an approximation depending on the current iterate, resulting in an approximated problem that can be solved more easily to yield the next iterate. By applying the semismooth Newton logic to solve~\eqref{eq:risk-averse_rootfinding}, given an iterate $v_k$, we form the following local approximation
\begin{equation}\label{eq:Newton_equation}
    r(v_k) + \xi(v_k,\,d)\,,
\end{equation}
where $\xi(v_k,\,d)$ belongs to a family of approximations $\Xi(v_k)$ where each element $\xi(v_k,\,\cdot)\in\Xi(v_k)$ is a function from $\mathbb{R}^n$ to $\mathbb{R}^n$ and should approximate ``well-enough'' $r(v_k + d) - r(v_k)$ around $d = 0$.  The other fundamental property of $\Xi(v_k)$ is that each $\xi(v_k,\,\cdot) \in \Xi(v_k)$ should lead to a system of equations that is easy to solve and, preferably, has a unique solution for $d$ so that the semismooth Newton iterate is well-defined. In particular, let $d_k$ be such that
\begin{equation}\label{eq:newton_approx_system}
     r(v_k) + \xi(v_k,\, d_k) = 0\,,
\end{equation}
the new Newton iterate is set to
\begin{equation}\label{eq:newton_iteration}
    v_{k+1}  = v_k + d_k\,,
\end{equation}
where $v_0\in \mathbb{R}^n$. The main difference with respect to the discussion in~\cite{gargiani_newton} is that we do not only consider linear approximations but~\eqref{eq:Newton_equation} could potentially be non-linear.
When it is possible to compute at least some of its elements, a common choice for $\Xi(v_k)$ is the set $\left\{ J _k d \quad | \quad J_k \in \partial r(v_k) \right\}$ where $\partial r(v_k)$ is the Clarke's generalized Jacobian of $r$ at $v_k$~\cite[Section 2.4]{izmailov_book}. With this choice, the resulting local model in~\eqref{eq:Newton_equation} is linear and, if $\partial r(v_k)$ is CD-regular~\cite[Remark 1.65]{izmailov_book}, it has a unique zero. 

Generally adopted requirements for selecting $\Xi(v)$ are that for all $v\in \mathbb{R}^n$ and for all $\xi(v,\,\cdot)\in \Xi(v)$
\begin{enumerate}
    \item $\xi(v,\,0) = 0$,
    \item $\xi(v,\,\cdot)$ is invertible and $d(v)= \xi^{-1}(v,\,-r(v))$ is the unique solution of $r(v) + \xi(v,\,d) = 0$.
\end{enumerate}
While (1) is a fundamental requirement for a Newton approximation, (2) could be further relaxed to account for scenarios where~\eqref{eq:Newton_equation} admits multiple solutions. These conditions are in general not sufficient to guarantee convergence of iteration~\eqref{eq:newton_iteration}. For that additional stronger conditions on $\Xi$ and structural properties inherent to $r$ should hold. In particular, in the following theorem we prove that, if $\Xi$ verifies certain structural properties in addition to the basic discussed requirements, then the sequence generated by a general semismooth Newton method enjoys a fast quadratic rate of convergence in a neighborhood of the solution.

\begin{thm}\label{th:quadratic_convergence_NA}
    Let $\left\{v_k \right\}$ be the sequence generated by a general semismooth Newton method and assume that for every $\xi(v,\,\cdot) \in \Xi(v)$ 
        \begin{enumerate}
            \item $\xi(v,\,\cdot)$ is a uniform Lipschitz homeomorphism for all $v\in\mathbb{R}^n$ in the sense of Definition 2.1.9 in~\cite{facchinei_book_vol1}, and $L_{\Xi}\geq 0$ is the Lipschitz modulus of the inverse,
            \item there exists a neighborhood 
 $\mathcal{U}_{v^*}'\subseteq\mathbb{R}^n$ and a constant $C_{\Xi}\geq 0$ such that for all $v\in\mathcal{U}_{v^*}'$
            \begin{equation}\label{eq:strong_semismoothness}
            \Vert r(v) + \xi(v,\,v^* - v) - r(v^*) \Vert \leq C_{\Xi}\,\Vert v - v^* \Vert^2\,.
            \end{equation}
        \end{enumerate}
        In addition, let
        \begin{equation}
            {\mathcal{U}''}_{v^*} = \left\{ v\in\mathcal{U}_{v^*}' \,\,\,\text{s.t.}\,\,\, \Vert v - v^* \Vert < 1/(L_{\Xi}C_{\Xi}) \right\}\,.
        \end{equation}
        Then, if $v_0 \in {\mathcal{U}''}_{v^*}$, 
        the sequence $\left\{ v_k \right\}$ is converging Q-quadratically to $v^*$.
\end{thm}
\noindent \textit{Proof.} Recall that $d_k = \xi^{-1}(v_k,\, -r(v_k))$ for any $k\geq 0$. Then
\begin{equation}\label{eq:upper_bounding}
    \begin{aligned}
        &\Vert v_{1}  -  v^* \Vert = \Vert v_{0} + \xi^{-1}(v_0,\, -r(v_0))  - v^* \Vert\\
        &= \Vert\xi^{-1}(v_0,\, -r(v_0)) - (v^* - v_0) \Vert\\
        &= \Vert\xi^{-1}(v_0,\, -r(v_0)) - \xi^{-1}(v_0,\,\xi(v_0,\,(v^* - v_0)) \Vert\\
        &\leq L_{\Xi}\, \Vert r(v_0) + \xi(v_0,\,v^* - v_0)\Vert\\
        &= L_{\Xi}\, \Vert r(v_0) + \xi(v_0,\,v^* - v_0) - r(v^*)\Vert\\
        &\overset{(a)}{\leq} L_{\Xi} C_{\Xi}\, \Vert v_0 - v^* \Vert^2\,,
    \end{aligned}
\end{equation}
where $(a)$ follows from the fact that, since $v_0\in{\mathcal{U}''}_{v^*}$, we can deploy relation~\eqref{eq:strong_semismoothness}.
It follows from~\eqref{eq:upper_bounding} that $\Vert v_1 - v^* \Vert\leq \delta_{\Xi} \Vert v_0 - v^*\Vert$ with $\delta_{\Xi} = L_{\Xi}C_{\Xi} \Vert v_0 - v^* \Vert$. Therefore also $v_1 \in \mathcal{U}'_{v^*}$ since $\delta_{\Xi}<1$. By repeating recursively the steps in~\eqref{eq:upper_bounding} and noticing with a similar reasoning that, for any $k>0$, $v_k \in \mathcal{U}'_{v^*}$, we obtain that 
\begin{equation}\label{eq:upperbound_k}
   \Vert v_{k+1} - v^* \Vert\leq L_{\Xi}C_{\Xi}\Vert v_k - v^* \Vert^2 \,.
\end{equation}
We can also show that all the following contraction factors are upper-bounded by $\delta_{\Xi}$, which results in the following upper-bound for any $k>0$
\begin{equation}
    \Vert v_k - v^* \Vert\leq \delta_{\Xi}^k\Vert v_0 - v^* \Vert\,.
\end{equation}
This implies that the sequence $\left\{ v_k \right\}$ with $v_0\in \mathcal{U}_{v^*}''$ is converging to $v^*$ with a rate that is at least $\delta_{\Xi}$. Going back to~\eqref{eq:upperbound_k} we can therefore conclude Q-quadratic convergence as
\begin{equation}
    \frac{\Vert v_{k+1} - v^* \Vert}{\Vert v_{k} - v^* \Vert^2} \leq L_{\Xi}C_{\Xi}\,.
\end{equation}~\qed

Similar results on local quadratic convergence of semismooth Newton methods using the concept of Newton approximation scheme and for general non-smooth root finding problems are known in the literature, \textit{e.g.}, see Section 7.2 in~\cite{facchinei_book}.

Condition $(2)$ in Theorem~\ref{th:quadratic_convergence_NA} is key to prove local quadratic convergence of a general semismooth Newton methods. In the literature a Newton approximation scheme which verifies this condition is also called \textit{strong Newton approximation scheme}~\cite[Definition 7.2.2]{facchinei_book}. Finally, notice that the results of Theorem~\ref{th:quadratic_convergence_NA} could be used to design \textit{ad-hoc} Newton approximation schemes that guarantee fast local rate of convergence.


We now introduce three different families of approximations for the risk-averse Bellman residual function~\eqref{eq:risk-averse_rootfinding}. These are later used in Sections~\ref{sec:SNM1}-\ref{sec:SNM3} to define different instances of the semismooth Newton method to solve~\eqref{eq:risk-averse_rootfinding}.
We start by introducing two auxiliary operators, which are then deployed for the definition of our approximations.

For $\tilde{v}\in\mathbb{R}^n$, we define $T^{\tilde{v}}:\mathbb{R}^n \rightarrow \mathbb{R}^n$ as 
\begin{equation}\label{eq:auxiliary_T}
    T^{\tilde{v}} v = \min_{\pi\in\Pi}\left\{ g_{\pi} + \gamma P^{\tilde{v}}_{\pi} v \right\}\,,
\end{equation}
where for every $\pi\in\Pi$ $P_{\pi}^{\tilde{v}}\in\mathbb{S}^{n\times n}$ is the row-stochastic matrix whose $i$-th row is defined as follows
\begin{equation}\label{eq:matrix_maximization}
    \left[P^{\tilde{v}}_{\pi}\right]_{i,:} \in \arg\max_{p\in \mathcal{P}_{\chi}(P(i,\pi(i)))} \langle p , \, \tilde{v} \rangle \quad \forall\, i\in\mathcal{S}\,.
\end{equation}
Notice that, by construction, we have that 
\begin{equation}\label{eq:relation_operators}
    T^{\tilde{v}} v \leq Dv \quad \forall \,v\in\mathbb{R}^n\,,
\end{equation}
since, by definition, $P_{\pi}^{\bar{v}} v\leq P_{\pi}^v v $ for any $\pi\in\Pi$.
In addition, $T^{\tilde{v}}\tilde{v} = D\tilde{v}$ by definition of the operators.
Similarly, for $\tilde{\pi}\in\Pi$, we introduce the operator $T^{\tilde{v}}_{\tilde{\pi}}:\mathbb{R}^n \rightarrow \mathbb{R}^n$ which maps $v\in\mathbb{R}^n$ to
\begin{equation}\label{eq:auxiliary_Tpi}
    T^{\tilde{v}}_{\tilde{\pi}}v = g_{\tilde{\pi}} + \gamma P_{\tilde{\pi}}^{\tilde{v}} v\,,
\end{equation}
where $P_{\tilde{\pi}}^{\tilde{v}}$ is defined as in~\eqref{eq:matrix_maximization}.

The operators in~\eqref{eq:auxiliary_T} and~\eqref{eq:auxiliary_Tpi} are interpretable as the risk-neutral Bellman operators for optimality and associated to the fixed policy $\tilde{\pi}$, respectively, for an MDP with~\eqref{eq:matrix_maximization} as transition probability map.
\begin{lem}\label{lm:auxiliary_operators}
 The operators in~\eqref{eq:auxiliary_Tpi} and~\eqref{eq:auxiliary_T} have~\eqref{eq:risk-neutral_DP_of} and the minimum of ~\eqref{eq:risk-neutral_DP_of} over $\pi\in \Pi$ for an MDP with~\eqref{eq:matrix_maximization} as transition probability map as unique fixed-points, respectively. They also verify the fundamental properties of monotonicity, constant shift and $\gamma$-contractivity in the infinity-norm~\cite[Chapter 1]{bertsekas_book}.  
\end{lem}
\noindent\textit{Proof.} The statements follow directly from the definition of the operators and the fundamental properties of the risk-neutral Bellman operators~\cite[Chapter 1]{bertsekas_book}.~$\square$

 We define the following three families of approximations for the risk-averse Bellman residual function:
\begin{enumerate}
    \item $\Xi_1$: $\xi_1(v, \,d) = d + Dv - T^v(v + d)$\,,
    \item $\Xi_2$: $\xi_2(v, \,d) = d + Dv - D_{\pi}(v+d)$ with ${D}_{\pi}v = {D}v$\,,
    \item $\Xi_3$: $\xi_3(v, \,d) = d + Dv - T_{\pi}^v(v + d)$ with ${D}_{\pi} v = {D} v$\,.
\end{enumerate}

\begin{lem}\label{lm:fundamental_properties_newton_approximations}
Let $v,\, d\in\mathbb{R}^n$ and define $\tilde{v} =  v + d$. All three approximations, $\Xi_1$, $\Xi_2$ and $\Xi_3$, verify the fundamental requirements for a Newton approximation. In addition 
\begin{enumerate}[label=(\roman*)]
    \item $r(v) + \xi_1(v,\,d) = \tilde{v} - {T}^v \tilde{v}$,
    \item $r(v) + \xi_2(v,\,d) = \tilde{v} - {D}_{\pi} \tilde{v}$,
    \item $r(v) + \xi_3(v,\,d) = \tilde{v} - {T}^v_{\pi}\tilde{v}$.
\end{enumerate}
\end{lem}
\textit{Proof.} Notice that, for a general $\Xi$, to prove that the second property holds it is sufficient to show that~\eqref{eq:Newton_equation} has a unique solution. We start by considering $\Xi_1$: 
\begin{enumerate}
    \item Let $\xi_1(v, \cdot) \in \Xi_1$. Then
    \begin{equation*}
    \xi_1(v, 0) = 0 + Dv - T^v(v +0 ) = Dv - T^vv = 0\,,
    \end{equation*}
    where the last equality follows from the fact that $T^v v = Dv$ by definition of the operators. 
    \item We compute explicitly the expression for the local model
    \begin{equation}\label{eq:derivation_newtonupdate1}
        \begin{aligned}
            r(v) + \xi_1(v,\,d) &= v - Dv + d + Dv - T^v(v+d)\\
            &= v + d - T^v(v+d) \\
            &= \tilde{v} - {T}^{v}\tilde{v}\,,
        \end{aligned}
    \end{equation}
where the last equality follows from the definition of $\tilde{v}$.
It follows directly from Lemma~\ref{lm:auxiliary_operators} that Equation~\eqref{eq:derivation_newtonupdate1} has a unique zero at the risk-neutral optimal cost for the perturbed MDP with transition probabilities computable as in~\eqref{eq:matrix_maximization}.
\end{enumerate}
We now consider $\xi_2(v,\,\cdot)\in\Xi_2$:
\begin{enumerate}
    \item By evaluating at $d=0$ we obtain
    \begin{equation*}
        \xi_2(v,\, 0) = 0 + Dv - D_{\pi}(v + 0) = Dv - D_{\pi} v = 0\,,
    \end{equation*}
    where the last equality follows from the fact that $\pi$ is a risk-averse greedy policy with respect to $v$.
    \item We compute the expression for the local model as follows
    \begin{equation}\label{eq:local_model2}
\begin{aligned}
    r(v) + \xi_2(v,\,d) &= v - {D}v + d  + Dv - D_{\pi}(v + d)\\
    &= v + d - D_{\pi}(v+d) \\
    &= \tilde{v} - {D}_{\pi}\tilde{v} \,,
\end{aligned}
\end{equation}
where the last equality follows from the definition of $\tilde{v}$. From Lemma~\ref{lemma:fundamental_properties_D} it follows that Equation~\eqref{eq:local_model2} has a unique zero at the risk-averse cost associated to policy $\pi$.
\end{enumerate}
Finally, for any $\xi_3(v,\,\cdot)\in\Xi_3(v)$:
\begin{enumerate}
    \item It can be proved by following analogous steps as for $(1)$ in the previous two cases and using the fact that $Dv = T_{\pi}^v v$ since $\pi$ is selected such that $D_{\pi}v = Dv$.
    \item We compute the local model as follows
    \begin{equation}\label{eq:newton_model3}
    \begin{aligned}
        r(v) + \xi_3(v,\,d) &= v - Dv + d + Dv - T_{\pi}^v(v + d) \\
        &= v + d - T_{\pi}^v(v + d) \\
        &= \tilde{v} - {T}^v_{\pi}\tilde{v}\,,
    \end{aligned}
    \end{equation}
    where the last equality follows from the definition of $\tilde{v}$.
    It follows directly from Lemma~\ref{lm:auxiliary_operators} that Equation~\eqref{eq:newton_model3} has a unique zero at the risk-neutral cost associated to policy $\pi$ for the perturbed MDP with transition probabilities computable as in~\eqref{eq:matrix_maximization}.
\end{enumerate}~\qed

\section{SNM I}\label{sec:SNM1}

In this section we consider $\Xi_1$ as family of approximations.
Let $v_k\in\mathbb{R}^n$ be the current iterate, by using the results of Lemma~\ref{lm:fundamental_properties_newton_approximations}, we obtain that
\begin{equation}
    \begin{aligned}
        r(v_k) &+ \xi_1(v_k,\,d_k) = v_{k+1} - {T}^{v_k}v_{k+1}\,.
    \end{aligned}
\end{equation}
Consequently, the choice of this local model results in the following Newton iteration
\begin{equation}\label{eq:newtonupdate_1}
    v_{k+1} = {T}^{v_k} v_{k+1}\,.
\end{equation}
Each cost-update therefore requires the solution of a risk-neutral MDP, where the transition matrices are computed according to~\eqref{eq:matrix_maximization}. For finite state-action MDPs we can solve~\eqref{eq:newtonupdate_1} in a finite number of iterations with risk-neutral policy iteration~\cite{bertsekas_book}. With this instance of the semismooth Newton method we solve the original risk-averse MDP by solving a sequence of risk-neutral MDPs with different transition probability matrices. We call this method \textit{Semismooth Newton Method I} (SNMI, see Algorithm~\ref{alg:Newton1}).

The following theorem shows that SNMI generates a non-decreasing sequence of iterates that converges globally to the minimizer of~\eqref{eq:obj_fun}.
\begin{thm}\label{th:snm1}
    The sequence $\left\{ v_k \right\}$ generated by the SNMI converges to the unique solution $v^*$ of~\eqref{eq:risk-averse_BE}. Moreover, the sequence is monotone, \textit{i.e.},
    \begin{equation}
        v_{k+1} \geq v_k \quad k=1,2,3,\dots\,.
    \end{equation}
\end{thm}
\noindent \textit{Proof.} Consider $k\geq 0$. Because of the monotonicity property and~\eqref{eq:relation_operators}, we obtain that
\begin{equation}
    ({T}^{v_k})^2 V \leq {T}^{v_{k}}{D}v \leq {D}^2 v\quad \forall\, v\in\mathbb{R}^n\,.
\end{equation}
Repeating this argument $W$-times leads to
\begin{equation}\label{eq:before_limit}
    ({T}^{v_k})^W v \leq {D}^W v \quad \forall\, V\in\mathbb{R}^n\,.
\end{equation}
Since both ${T}^{v_k}$ and ${D}$ have unique fixed-points (by Lemma~\ref{lm:auxiliary_operators} and Lemma~\ref{lemma:fundamental_properties_D}), we can take the limit for $W\rightarrow \infty$ of~\eqref{eq:before_limit} and we obtain that
\begin{equation}\label{eq:sequence_boundedabove}
    v_{k+1} \leq v^*\,,
\end{equation}
since, by definition, $v_{k+1}$ is the fixed-point of $T^{v_k}$ (see Equation~\eqref{eq:newtonupdate_1}).
We can therefore conclude that the sequence $\left\{ v_k \right\}$ is bounded above by $v^*$. In addition, for any $k>0$ we have
\begin{equation}\label{eq:sequence_nondecreasing}
\begin{aligned}
    v_k &= {T}^{v_{k-1}} v_k \leq {D} v_k \\ 
    &= {T}^{v_k}v_k \overset{(a)}{\leq} ({T}^{v_k})^W v_k \overset{W\rightarrow \infty}{\longrightarrow} v_{k+1}\,,
\end{aligned}
\end{equation}
where $(a)$ follows from the monotonicity of the ${T}^{V_k}$-operator and the limit by contractivity (Lemma~\ref{lm:auxiliary_operators}) and the fact that $v_{k+1}$ is the fixed point of $T^{v_k}$ by definition.
Since the sequence is non-decreasing~\eqref{eq:sequence_nondecreasing} and bounded above~\eqref{eq:sequence_boundedabove}, it admits a limit $\hat{v}\in\mathbb{R}^n$. Taking the limit for $k\rightarrow \infty$ in~\eqref{eq:sequence_nondecreasing} and since, by contractivity of ${D}$ (see Lemma~\ref{lemma:fundamental_properties_D}), $\lim_{k\rightarrow \infty} {D}v_k = {D}\lim_{k\rightarrow \infty} v_k$, we obtain that
\begin{equation}
    \hat{v} \leq {D}\hat{v} \leq \hat{v} \,.
\end{equation}
Since ${D}$ has $v^*$ as unique fixed-point, then we can conclude that the sequence $\left\{ v_k \right\}$ is also converging to $v^*$.~\qed

Finally, note that for SNMI the condition~\eqref{eq:strong_semismoothness} in Theorem~\ref{th:quadratic_convergence_NA} reduces to requiring that there exists $\mathcal{U}'_{v^*, 1}$ and $C_1\geq 0$ such that
\begin{equation}\label{eq:ss_snm1}
    \Vert v^* - T^v v^*\Vert \leq C_1\Vert v - v^* \Vert^2 \quad \forall\, v \in \mathcal{U}'_{v^*, 1}\,.
\end{equation}
This condition trivially holds when the operator is piecewise affine since piecewise affine functions are strongly semismooth~\cite[Proposition 7.4.7]{facchinei_book}.
\begin{algorithm}
\caption{Semismooth Newton Method I}
\label{alg:Newton1}
\begin{algorithmic}[1]
\Require $v_0 \in \mathbb{R}^n$
\State $k\leftarrow 0$
\While{convergence criterion is not reached}
\For {$\pi \in \Pi$}
\State compute $P_{\pi}^{v_k}$ using~\eqref{eq:matrix_maximization}
\EndFor
\State solve $v_{k+1} = {T}^{v_k} v_{k+1} $
\State $k \leftarrow k+1$
\EndWhile
\end{algorithmic}
\end{algorithm}

\section{SNM II}\label{sec:SNM2}

If we consider the second class of functions $\Xi_2$, by using the results of Lemma~\ref{lm:fundamental_properties_newton_approximations}, we obtain the following local model at iteration $k\geq 0$
\begin{equation}
\begin{aligned}
    r(v_k) &+ \xi_2(v_k,\,d_k) = v_{k+1} - {D}_{\pi_{k+1}}v_{k+1} \,,
\end{aligned}
\end{equation}
where ${D}_{\pi_{k+1}}v_k = {D}v_k$.
This choice of model leads to the following Newton iteration
\begin{equation}\label{eq:newton_update2}
    v_{k+1} = {D}_{\pi_{k+1}} v_{k+1}\,.
\end{equation}
Because of the potential non-linearity and non-smoothness of operator~\eqref{eq:risk-averse_BE_pi}, ~\eqref{eq:newton_update2} requires the solution of a non-trivial non-smooth and non-linear system of equations whose structure depends on the specific choice of the risk-transition map. We call this method \textit{Semismooth Newton Method II} (SNMII, see Algorithm~\ref{alg:snm2}).

We can look at~\eqref{eq:newton_update2} as the evaluation of the cost of policy $\pi_{k+1}$ for the considered risk-averse MDP. The latter is a risk-averse greedy policy with respect to the current iterate $v_k$. Therefore, this instance of semismooth Newton mimics in the risk-averse scenario the risk-neutral policy iteration method, as it alternates a policy improvement step and a policy evaluation step. We note that iteration~\eqref{eq:newton_update2} is identical to the iteration proposed by Ruszczyński in~\cite[Section 8]{risk_averse_DP}. In particular, iterations (46) and (47) in~\cite{risk_averse_DP} can be obtained by considering $\Xi_2$ as family of local approximations in~\eqref{eq:Newton_equation}.   

Unlike the standard policy iteration scheme, in our setting the policy evaluation step requires the solution of a non-linear and non-smooth system of equations. Ruszczyński proposes the deployment of a semismooth Newton method for this task giving rise to the inner iteration
\begin{equation}\label{eq:newton_inneriteration}
    \tilde{v}_{i+1} = g_{\pi_{k+1}} + \gamma P_{\pi_{k+1}}^{\tilde{v}_i} \tilde{v}_{i+1}\,,
\end{equation}
where $P_{\pi_{k+1}}^{\tilde{v}_i}$ is defined as in~\eqref{eq:matrix_maximization} and $\tilde{v}_0\in\mathbb{R}^n$~\cite[Section 9]{risk_averse_DP}. It can be shown that~\eqref{eq:newton_inneriteration} converges globally to the unique solution  of~\eqref{eq:newton_update2}~\cite[Theorem 7]{risk_averse_DP}. Global convergence of the sequence $\left\{v_k \right\}$ generated by Algorithm~\ref{alg:snm2} to the fixed-point of~\eqref{eq:risk-averse_BE} then follows from the policy improvement property~\cite[Theorem 7]{risk_averse_DP}; this, however, is guaranteed only if~\eqref{eq:newton_update2} is solved exactly at every iteration, a condition that would unrealistically require an infinite number of iterations of~\eqref{eq:newton_inneriteration}. In practice, the inner iterations are terminated when a certain tolerance level is reached.
\begin{algorithm}
\caption{Semismooth Newton Method II}
\label{alg:snm2}
\begin{algorithmic}
\Require $v_0 \in \mathbb{R}^n$
\State $k\leftarrow 0$
\While{convergence criterion is not reached}
\For {$\pi \in \Pi$}
\State compute $P_{\pi}^{v_k}$ using~\eqref{eq:matrix_maximization}
\EndFor
\State extract $\pi_{k+1} $ such that ${D}_{\pi_{k+1}} v_k = {D} v_k $ 
\State solve $v_{k+1} = {D}_{\pi_{k+1}} v_{k+1} $
\State $k \leftarrow k+1$
\EndWhile
\end{algorithmic}
\end{algorithm}

Exploring the parallel to policy iteration further, we could mimic optimistic policy iteration (OPI)~\cite[Section 2.3]{bertsekas_book}, and, instead of solving~\eqref{eq:newton_update2} exactly, produce an estimate of $v^{\pi_{k+1}}$ by performing a finite number of risk-averse value iteration steps with the operator ${D}_{\pi_{k+1}}$ starting from the current iterate $v_k$. This would lead to the iteration
\begin{equation}
    v_{k+1} = \left({D}_{\pi_{k+1}}\right)^{w_k} v_k\,,
\end{equation}
where $\pi_{k+1}$ is the risk-averse greedy policy with respect to $v_k$ and $w_k > 0$ is the number of risk-averse value iteration steps used (Algorithm~\ref{alg:riskaverse_opi}). 
The following theorem characterizes the convergence properties of the resulting sequence $\left\{ v_k \right\}$.
\begin{thm}\label{th:ra_opi}
    The sequence $\left\{ v_k \right\}$ generated by Algorithm~\ref{alg:riskaverse_opi} converges to $v^*$ for any $v_0\in \mathbb{R}^n$.  
\end{thm}
\noindent \textit{Proof.} We start by selecting a scalar $b$ such that ${D}\bar{v}_0 \leq \bar{v}_0$ with $\bar{v}_0 = v_0 + be$. Notice that the $b$-shift does not affect the generated sequence of policies. Indeed, thanks to Lemma~\ref{lm:constant_shift_D}, it is possible to verify that if we replace $v_0$ with $\bar{v}_0$ Algorithm~\ref{alg:riskaverse_opi} generates the same sequence of policies $\left\{ \pi_{k+1}\right\}$. We now show convergence of $\left\{\bar{v}_k \right\}$ to $v^*$ and then use this result to prove the convergence of the original sequence $\left\{ v_k\right\}$ to $v^*$.

We start the proof by recalling that ${D}_{\pi_1} \bar{v}_0 = {D}\bar{v}_0 \leq \bar{v}_0$, where the equality follows from the fact that $\pi_1$ is the risk-averse greedy policy associated to $\bar{v}_0$. From this relation together with the monotonicity of the risk-averse operators, we obtain that for $w=1,2,\dots$
\begin{multicols}{2}
  \begin{equation*}
    \left({D}\right)^w \bar{v}_0 \leq \left({D}\right)^{w-1} \bar{v}_0
  \end{equation*}
  \break
  \begin{equation*}
    \left({D}_{\pi_1}\right)^w \bar{v}_0 \leq \left({D}_{\pi_1}\right)^{w-1} \bar{v}_0
  \end{equation*}
\end{multicols}

By the definition of the risk-averse greedy policy, this leads to
\begin{equation}
    \begin{aligned}
        {D}_{\pi_2} \bar{v}_1 &= {D} \bar{v}_1\leq {D}_{\pi_1} \bar{v}_1
        \overset{(a)}{=} \left({D}_{\pi_1}\right)^{w_0 + 1} \bar{v}_0\\
        & \overset{(b)}{\leq} \left({D}_{\pi_1}\right)^{w_0} \bar{v}_0 = \bar{v}_1 \overset{(c)}{\leq} {D}_{\pi_1} \bar{v}_0 = {D} \bar{v}_0\,,
    \end{aligned}
\end{equation}
where $(a)$ follows from the fact that $\bar{v}_{k+1} = \left({D}_{\pi_{k+1}}\right)^{w_k}\bar{v}_k$ for any $k\geq 0$.
This argument can be continued to show that $\bar{v}_{k+1} \leq {D}\bar{v}_{k}$ for $k\geq 0$ and therefore that 
\begin{equation}\label{eq:ineq_1}
\bar{v}_{k+1} \leq \left({D}\right)^{k+1}\bar{v}_{0}\,.
\end{equation}
Recall that for a general $k\geq 0$
\begin{equation}\label{eq:vk_uwrapped}
    \bar{v}_{k+1} = \left({D}_{\pi_{k+1}}\right)^{w_k}\left({D}_{\pi_{k}}\right)^{w_{k-1}}\dots\left({D}_{\pi_1}\right)^{w_0}\bar{v}_0\,.
\end{equation} 
Starting from ${D}_{\pi_1} \bar{v}_0 = {D}\bar{v}_0$ and applying the ${D}_{\pi_1}$-operator on both sides leads to
\begin{equation}\label{eq:step1}
    \left({D}_{\pi_1}\right)^2 \bar{v}_0 = {D}_{\pi_1}{D}\bar{v}_0 \overset{(a)}{\geq} \left({D}\right)^2 \bar{v}_0\,,
\end{equation}
where $(a)$ follows from the fact that for any $V\in \mathbb{R}^n$ and for any $\pi\in\Pi$
\begin{equation}\label{eq:versus_operators}
    {D}_{\pi}v \geq {D}v\,.
\end{equation}
By repeating this argument $w_0$-times we obtain that
\begin{equation}\label{eq:ineq_tmp}
    \left({D}_{\pi_1}\right)^{w_0} \bar{v}_0 \geq \left({D}\right)^{w_0} \bar{v}_0\,.
\end{equation}
Recall that $\bar{v}_1 = \left({D}_{\pi_1}\right)^{w_0} \bar{v}_0$. If we set $\hat{v}_1 = \left({D}\right)^{w_0} \bar{V}_0$, we can rewrite~\eqref{eq:ineq_tmp} as $\bar{v}_1\geq \hat{v}_1$. By applying the ${D}_{\pi_2}$-operator on both sides of this inequality, we obtain that
\begin{equation}\label{eq:ineq_cascate}
    {D}_{\pi_2} \bar{v}_1 \overset{(a)}{\geq} {D}_{\pi_2} \hat{v}_1 \overset{(b)}{\geq} {D}\hat{v}_1\,,
\end{equation}
where $(a)$ follows from the monotonicity of the ${D}_{\pi_2}$-operator which makes its application inequality invariant, and $(b)$ from~\eqref{eq:versus_operators}. 
Continuing with this same argument for $w_1$-times we can conclude that
\begin{equation}\label{eq:finalstep}
    \left({D}_{\pi_2}\right)^{w_1} \bar{v}_1 \geq \left({D}\right)^{w_1}\hat{v}_1\,.
\end{equation}
We can now set $\hat{v}_2 = \left({D}\right)^{w_1}\hat{v}_1$ and rewrite~\eqref{eq:finalstep} as $\bar{v}_2\geq \hat{v}_2$.
Starting from~\eqref{eq:finalstep} and repeating $(k-2)$-times the same steps that have led from~\eqref{eq:ineq_tmp} to~\eqref{eq:finalstep}, we can conclude that
\begin{equation}\label{eq:ineq_2}
\begin{aligned}
 \bar{v}_{k+1} &\overset{(a)}{=}  \left({D}_{\pi_{k+1}}\right)^{w_k}\dots \left({D}_{\pi_1}\right)^{w_0}\bar{v}_0 \\
 &\geq \hat{v}_{k+1} \overset{(b)}{=} \left({D}\right)^{w_0 + \dots + w_k} \bar{v}_0\,,
\end{aligned}
\end{equation}
where $(a)$ follows from~\eqref{eq:vk_uwrapped} and $(b)$ from the fact that $\hat{v}_1 = \left({D}\right)^{w_0}\bar{v}_0$ and $\hat{v}_{k+1} = \left({D}\right)^{w_k}\hat{v}_k$ for $k>0$.
By combining~\eqref{eq:ineq_1} and~\eqref{eq:ineq_2}, we obtain that for any $k>0$
\begin{equation}\label{eq:final_ineq}
    \left({D}\right)^{w_0 + \dots + w_{k-1}} \bar{v}_0 \leq \bar{v}_{k} \leq \left({D}\right)^{k}\bar{v}_0\,.
\end{equation}
Recall that $v^*$ is the unique fixed point of ${D}$ which is a contractive operator.
By taking the limit of~\eqref{eq:final_ineq} for $k\rightarrow \infty$, we can therefore conclude that
\begin{equation}\label{eq:final_result}
    \lim_{k\rightarrow \infty} \bar{v}_{k} = v^*\,,
\end{equation}
Since $\lim_{k\rightarrow \infty} \bar{v}_{k} - v_k = 0$ (for details see the proof of Proposition 9 in~\cite{gargiani_ipi}), from~\eqref{eq:final_result} we can also conclude that $\left\{ v_k \right\}$ converges to $v^*$.~\qed

\begin{algorithm}
\caption{Risk-Averse Optimistic Policy Iteration}
\label{alg:riskaverse_opi}
\begin{algorithmic}
\Require $v_0 \in \mathbb{R}^n$
\State $k\leftarrow 0$
\While{convergence criterion is not reached}
\For {$\pi \in \Pi$}
\State compute $P_{\pi}^{v_k}$ using~\eqref{eq:matrix_maximization}
\EndFor
\State extract $\pi_{k+1}$ such that ${D}_{\pi_{k+1}} v_k = {D} v_k $ 
\State $v_{k+1} = \left({D}_{\pi_{k+1}}\right)^{w_k} v_{k} $
\State $k \leftarrow k+1$
\EndWhile
\end{algorithmic}
\end{algorithm}
Notice that, as in the risk-neutral scenario, if we set $w_k= 1$ for all $k\geq 0$, then we recover risk-averse VI, while if we hypothetically run an infinite number of iterations and therefore solve exactly~\eqref{eq:newton_update2} at each iteration, we would recover Algorithm~\ref{alg:snm2}.

Finally, note that for SNMII the condition~\eqref{eq:strong_semismoothness} in Theorem~\ref{th:quadratic_convergence_NA} reduces to requiring that there exists $\mathcal{U}'_{V^*, 2}$ and $C_2\geq 0$ such that
\begin{equation}\label{eq:ss_snm2}
    \Vert v^* - D_{\pi} v^*\Vert \leq C_2\Vert v - v^* \Vert^2 \quad \forall\, v \in \mathcal{U}'_{v^*, 2}\,.
\end{equation}
This condition trivially holds when the operator is piecewise affine since piecewise affine functions are strongly semismooth~\cite[Proposition 7.4.7]{facchinei_book}.

\section{SNM III}\label{sec:SNM3}

Finally, we deploy the third class of functions leading to
\begin{equation}\label{eq:newtonupdate_3}
    v_{k+1} = {T}_{\pi_{k+1}}^{v_k} v_{k+1}\,, 
\end{equation}
where $\pi_{k+1}$ is the risk-averse greedy policy with respect to $v_k$ satisfying ${D}_{\pi_{k+1}} v_k = {D} v_k$. Note that while~\eqref{eq:newtonupdate_1} and~\eqref{eq:newton_update2} require the solution of a non-linear and non-smooth system of equations,~\eqref{eq:newtonupdate_3} is linear. This in general greatly facilitates its solution, as the numerical analysis literature provides many efficient direct and indirect methods to solve linear systems with non-singular coefficient matrices~\cite{Hackbusch94, saad}. Note also that the update in~\eqref{eq:newtonupdate_3} corresponds to the evaluation of the cost associated to the policy $\pi_{k+1}$ for a risk-neutral MDP with the transition matrix $P_{\pi_{k+1}}^{v_k}$ defined in~\eqref{eq:matrix_maximization}. We call this method \textit{Semismooth Newton Method III} (SNMIII, see Algorithm~\ref{alg:Newton3}).
\begin{algorithm}
\caption{Semismooth Newton Method III}
\label{alg:Newton3}
\begin{algorithmic}
\Require $v_0 \in \mathbb{R}^n$
\State $k\leftarrow 0$
\While{convergence criterion is not reached}
\For {$\pi \in \Pi$}
\State compute $P_{\pi}^{v_k}$ using~\eqref{eq:matrix_maximization}
\EndFor
\State extract $\pi_{k+1}$ such that ${D}_{\pi_{k+1}} v_k = {D} v_k$ 
\State solve $v_{k+1} = {T}_{\pi_{k+1}}^{v_k} v_{k+1} $
\State $k \leftarrow k+1$
\EndWhile
\end{algorithmic}
\end{algorithm}
As for the risk-neutral policy iteration method, SNMIII attempts to find the root of~\eqref{eq:risk-averse_BE} by solving a sequence of linear systems, each corresponding to the evaluation of a policy for an MDP with a perturbed transition map. The latter is recomputed at every iteration based on the current cost-vector. This method can also be interpreted as an inexact variant of SNMII, where only one iteration of~\eqref{eq:newton_inneriteration} with warm-starting is deployed to approximately solve~\eqref{eq:newton_update2} at every iteration. 

Finally, note that for SNMIII the condition~\eqref{eq:strong_semismoothness} in Theorem~\ref{th:quadratic_convergence_NA} reduces to requiring that there exists $\mathcal{U}'_{V^*, 3}$ and $C_3\geq 0$ such that
\begin{equation}\label{eq:ss_snm3}
    \Vert V^* - T_{\pi}^v v^*\Vert \leq C_3\Vert v - v^* \Vert^2 \quad \forall\, v \in \mathcal{U}'_{v^*, 3}\,.
\end{equation}
This condition trivially holds when the operator is piecewise affine since piecewise affine functions are strongly semismooth~\cite[Proposition 7.4.7]{facchinei_book}.

To establish convergence, the arguments of Theorem~\ref{th:snm1} and Theorem~\ref{th:ra_opi} are not applicable as, without additional assumptions, the local linearization breaks the monotonicity structure. While global convergence remains an open question for SNMIII, in the following theorem we characterize its local convergence properties for the class of piecewise continuously differentiable risk-measures by establishing a connection with piecewise smooth Newton~\cite[Chapter 7]{facchinei_book}. We refer to~\cite[Definition 1]{gargiani_newton} for a formal definition of piecewise continuously differentiable functions of order 1 (PC$^1$) and active pieces of a PC$^1$ function. 

\begin{thm}\label{th:Qsuperlinear_contraction}
For every $v_0 \in \mathbb{R}^n$ there exists a unique sequence $\left\{v_k \right\}_{k=0}^{\infty}$ satisfying the iteration of Algorithm~\ref{alg:Newton3}. Moreover, if $\chi :\mathbb{R}^n \times \mathbb{S}^n \rightarrow \mathbb{R}$ is piecewise continuously differentiable of order 1 in the first argument, then there exists a neighborhood $\mathcal{U}_{v^*}\subseteq\mathbb{R}^n$ of $v^*$ such that, for all $v_0\in \mathcal{U}_{v^*}$, the sequence converges Q-superlinearly to $v^*$.
\end{thm}
\textit{Proof}. Since in the considered setting $\Pi$ is a finite set and, by assumption, $\chi$ is PC$^1$, then we can conclude that the risk-averse Bellman residual function is also PC$^1$. For any $\tilde{v}\in \mathbb{R}^n$, the active pieces at $\tilde{v}$ are given by the following set of functions 
\begin{equation}
    \mathcal{Z}(\tilde{v}) = \left\{ v - {T}_{\pi}^{\tilde{v}}v \,\,\,\,\text{with}\,\,\,\, {D}_{\pi}\tilde{v} = {D}\tilde{v},\,v\in\mathbb{R}^n\right\}\,.
\end{equation}
The piecewise smooth Newton method~\cite[Chapter 7]{facchinei_book} considers the following set as approximation scheme for~\eqref{eq:risk-averse_rootfinding}
\begin{equation}
    \Xi(\tilde{v}) = \left\{ \left(I - \gamma P_{\pi}^{\tilde{v}}\right)d\quad \text{with}\quad {D}_{\pi}\tilde{v} = {D}\tilde{v} \right\}\,,
\end{equation}
which are the linear functions obtained by considering the Jacobians of the active pieces as coefficient matrices. The resulting local model at iteration $k\geq 0$ is 
\begin{equation}\label{eq:local_model_PCNM}
\begin{aligned}
    &r(v_k) + \left(I - \gamma P_{\pi_{k+1}}^{v_k} \right)  d_k
    \\&= v_k - {D}v_k + \left(I - \gamma P_{\pi_{k+1}}^{v_k} \right)  (v_{k+1} - v_k)
    \\&\overset{(a)}{=} v_{k+1} - g_{\pi_k} - \gamma P_{\pi_{k+1}}^{v_k}v_{k+1}
    \,,
\end{aligned}
\end{equation}
where $\pi_{k+1}$ is such that ${D}_{\pi_{k+1}}v_k = {D}v_k$ and $(a)$ follows from the fact that ${D}v_k = g_{\pi_{k+1}} + \gamma P_{\pi_{k+1}}^{v_k} v_k$ by definition of the ${D}$-operator.
Equation~\eqref{eq:local_model_PCNM} gives rise to the following iteration
\begin{equation}
    v_{k+1} = {T}_{\pi_{k+1}}^{v_k} v_{k+1}\,,
\end{equation}
which proves that SNMIII is an instance of the piecewise smooth Newton method for~\eqref{eq:risk-averse_rootfinding}. Since any matrix of the form $I - \gamma P$ with $P\in \mathbb{S}^{n\times n}$ is non-singular~\cite[Lemma 1]{gargiani_igmres}, by using the results of Theorem 7.2.15 in~\cite[Chapter 7]{facchinei_book} we can conclude that the sequence $\left\{ v_k \right\}$ generated by Algorithm~\ref{alg:Newton3} is unique and locally Q-superlinearly converging to $v^*$.~\qed

The main additional assumption of Theorem~\ref{th:Qsuperlinear_contraction} is verified for instance when CVaR is considered as risk-measure, and we show that formally in the following Proposition.
\begin{prop}\label{prop:cvar_pa}
  CVaR is a piecewise affine function in its first-argument.
\end{prop}
\noindent \textit{Proof.} When considering CVaR, the RHS in~\eqref{eq:dual_rewrite} reduces to a linear program where the feasible set is a polytope and it is independent from $V$ (see~\cite{NIPS2015_64223ccf}). 
Consequently, since there will be always a vertex of the polytope achieving the maximum, we could reduce the arg-maximization to the finite set comprising the vertices of the feasible region. For any $(s, a) \in \mathcal{S}\times \mathcal{A}$ and $V\in \mathbb{R}^n$, we can therefore re-write~\eqref{eq:dual_rewrite} as
\begin{equation}
    \chi(V, P(s, a)) = \max_{\tilde{p} \in \mathcal{V}(s, a)} \langle V, \tilde{p} \rangle\,, 
\end{equation}
where $\mathcal{V}(s, a) = \left\{ \tilde{p}_1,\dots, \tilde{p}_{\vert \mathcal{V} \vert} \right\}$ is the set comprising the vertices of the polytope defining the feasible region.
Therefore, since the set of admissible policies is finite, for any fixed $(s, a)$-pair, we can re-write $\chi$ as the maximum of a finite-collection of affine functions, which is a piecewise affine function by definition~\cite[Definition 4.1.3]{facchinei_book_vol1}.~\qed

To conclude, Proposition~\ref{prop:cvar_pa} allows one to conclude that conditions~\eqref{eq:ss_snm1},~\eqref{eq:ss_snm2} and~\eqref{eq:ss_snm3} hold when CVaR is deployed as risk-measure in~\eqref{eq:risk-averse_BE} and, therefore, we can deploy the results of Theorem~\ref{th:quadratic_convergence_NA} to conclude that the proposed SN methods enjoy local Q-quadratic contraction.

\section{Computational Considerations}

At each iteration of Algorithms~\ref{alg:Newton1}-\ref{alg:Newton3} $\vert \mathcal{S} \vert \times \vert \mathcal{A} \vert$ convex programs need to be solved for the evaluation of the operators which define the iteration of the methods. While the objective function of these programs is always linear, the feasible region is always convex, but its specific structure depends on the choice of the risk-measure. Beside these convex programs, each method requires a different amount of computation per iteration. These convex programs are independent of each other and can be solved in parallel. In particular, SNMI requires the solution of a risk-neutral MDP, which can be executed in less than $m^n$ iterations with risk-neutral policy iteration. With SNMII every iteration of the inner semismooth Newton method requires the solution of a linear system of dimension $\vert \mathcal{S} \vert$ and the re-evaluation of the transition map for the fixed policy which requires the solution of $\vert \mathcal{S} \vert$ convex programs. To enjoy the mathematical guarantees of global convergence theoretically an infinite number of iterations are required. In practice the number of iterations is truncated to a finite number, for example, using a parametric stopping condition. Finally, SNMIII only requires the solution of a linear system of dimension $\vert \mathcal{S} \vert$ per iteration.

As final remark on the computational consideration, it is worth to stress that the solution of the convex programs for the computation of the perturbed transition map in~\eqref{eq:matrix_maximization} is fully parallelizable as the programs are completely disjoint. 

\section{Numerical Evaluation}\label{sec:numerical-experiments}

To evaluate the numerical performance of our algorithms  we consider Conditional Value at Risk (CVaR) as the risk-measure, see, for example,~\cite{Rockafellar2000OptimizationOC}, \cite{Pflug2000}, \cite{sarykalin}. 
 CVaR  is a popular mathematical tool for managing risk, among others, in finance for portfolio optimization and asset liability management~\cite{presentation_CVaR}. Unlike more classical risk-measures such as the Value at Risk (VaR)~\cite{var_paper}, CVaR is a coherent risk-measure, which makes the optimization more stable and efficient~\cite{sarykalin} and allows us to deploy our framework. Our goal however is not that of comparing different risk measures, but rather that of empirically studying the performance of the discussed algorithms in terms of numerical stability and convergence. 
 
 There are multiple equivalent definitions of CVaR, we work with the one that characterizes it via the following optimization problem
\begin{equation*}
    \text{CVaR}(v(s_{t+1})) \! = \!\min_{z\in\mathbb{R}}\left\{ \! z + \frac{1}{\zeta}\mathbb{E}\left[ \left(v(s_{t+1}) - z\right)_{+} \,|\,\mathcal{F}_t \right] \! \right\} ,
\end{equation*}
where $\zeta\in(0,\,1)$ is the confidence-level parameter and $\mathcal{F}_t$ is the collection of the history.
The associated risk-transition map is
\begin{equation}\label{eq:cvar_risktransitionmap}
    \chi(v,\,P(s_t,\,a_t)) = \min_{z\in\mathbb{R}} \left\{ z + \frac{1}{\zeta} \langle (v - ze)_{+}, \, P(s_t,\,a_t) \rangle \right\} \,.
\end{equation}
The dual representation of~\eqref{eq:cvar_risktransitionmap} is as in~\eqref{eq:dual_rewrite} where
\begin{equation}
\begin{aligned}\label{eq:perturbed_distribution_set}
    \mathcal{P}_{\chi}(P(s,a))= \left\{ \tilde{p} \in\mathbb{S}^n\,|\,\tilde{p}= P(s,a) \odot \xi,\,\xi\in  [1,\,1/\zeta]^n \right\}.
\end{aligned}
\end{equation}
For completeness, we report the linear program associated with $v\in\mathbb{R}^n$ and $(s,a)\in\mathcal{S}\times\mathcal{A}$:
\begin{equation}
\begin{aligned}
    \left[P_{a}^{{v}}\right]_{s,:} = \arg&\max_{\tilde{p}\in\mathbb{R}^n} \,\,\langle v,\,\tilde{p} \rangle \\
    &\!\! \text{s.t.} \,\,\,\,\,\,\,\,\, \langle\tilde{p},\,e \rangle =1\\
    &\!\! \phantom{\text{s.t.}} \,\,\,\,\,\,\,\,\,
    0\leq \tilde{p} \leq P(s,a)/\zeta\,,
\end{aligned}    
\end{equation}
where the inequalities are interpreted component-wise.

Following~\cite{NIPS2015_64223ccf}, the CVaR of the random variable $v(s_{t+1})$ can be interpreted as the worst-case expectation of $v(s_{t+1})$ under a perturbed distribution $\tilde{p}(s_t,a_t)=\langle P(s_t,a_t),\xi \rangle$ where the perturbation $\xi$ depends on the confidence level $\zeta$. The perturbation could be considered as a modeling error, giving CVaR a robustness interpretation~\cite{NIPS2015_64223ccf, risk_averse_DP}.

Note that since~\eqref{eq:perturbed_distribution_set} is a polytope, the computation of~\eqref{eq:matrix_maximization} for a given policy $\pi\in\Pi$ reduces to the solution of $\vert \mathcal{S} \vert$ linear programs. Moreover, the optimization over policies needed by all algorithms is equivalent to the optimization over actions, for a total of at most $\vert \mathcal{S} \vert \times \vert \mathcal{A} \vert$ linear programs. 

As convergence criterion we used the infinity norm of the risk-averse Bellman residual function: if this quantity falls below a certain tolerance level then the algorithm returns the current iterate $v_k$ and its associated greedy policy as estimates of $v^*$ and $\pi^*$, respectively. An analogous criterion is deployed as stopping condition of the inner loop in SNMI and SNMII. In our benchmarks we used $10^{-6}$ for both the tolerance levels. The implementation is done in Python and we use cvxopt to solve the linear programs~\cite{cvxopt}. All benchmarks are run on Intel(R) Core(TM) i7-10750H CPU @ 2.60GHz.
The code is publicly available at \url{https://github.com/gmatilde/riskaverse\_DP.git}.

\begin{figure}
\centering
\begin{subfigure}{0.45\textwidth}
    \centering
    \includegraphics[width=\linewidth]{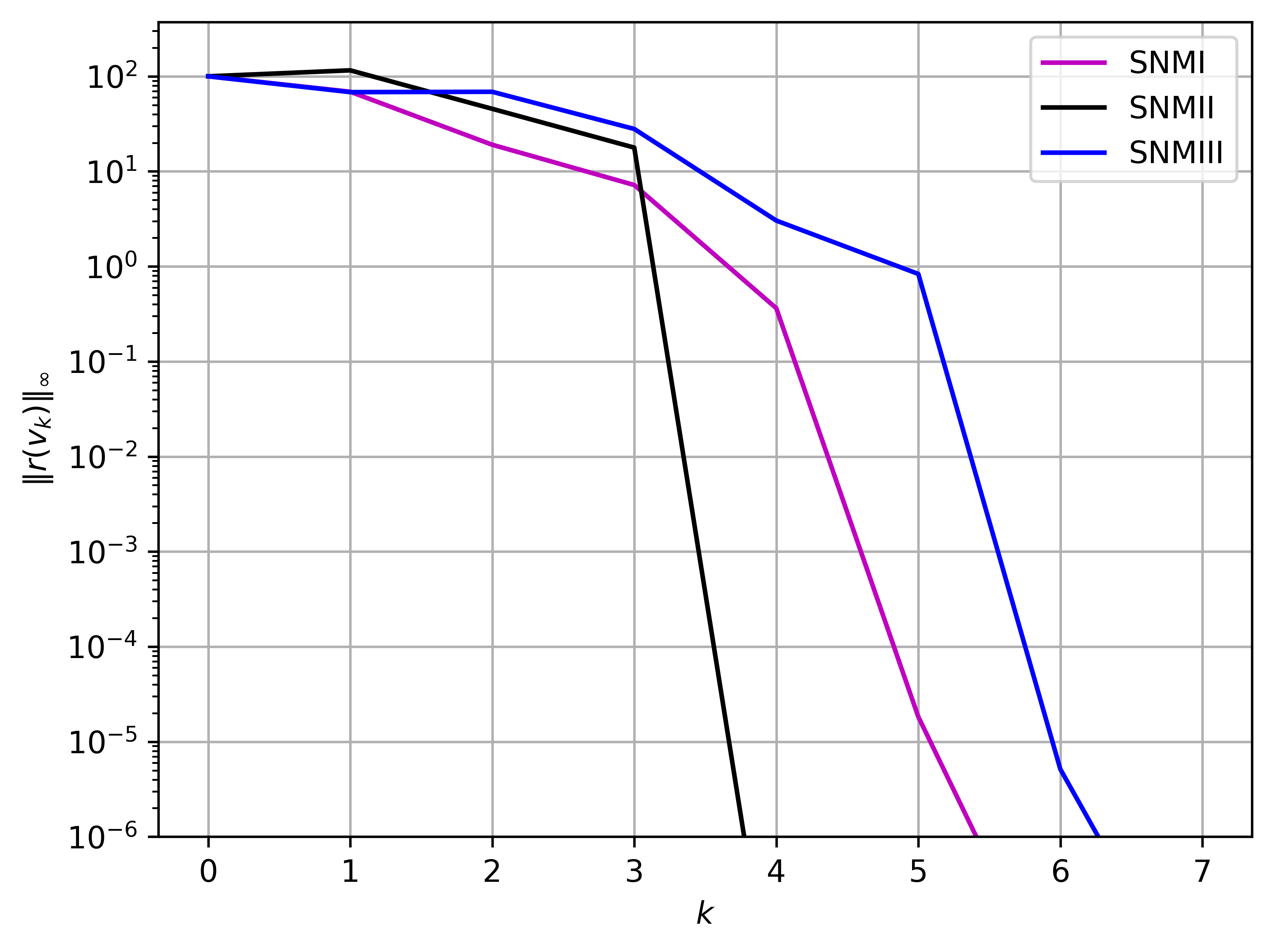}
    \caption{Infinity norm of the Bellman residual function versus the number of iterations for SNM methods.}
    \label{fig:SNM}
\end{subfigure}
\begin{subfigure}{0.45\textwidth}
    \centering
    \includegraphics[width=\linewidth]{
   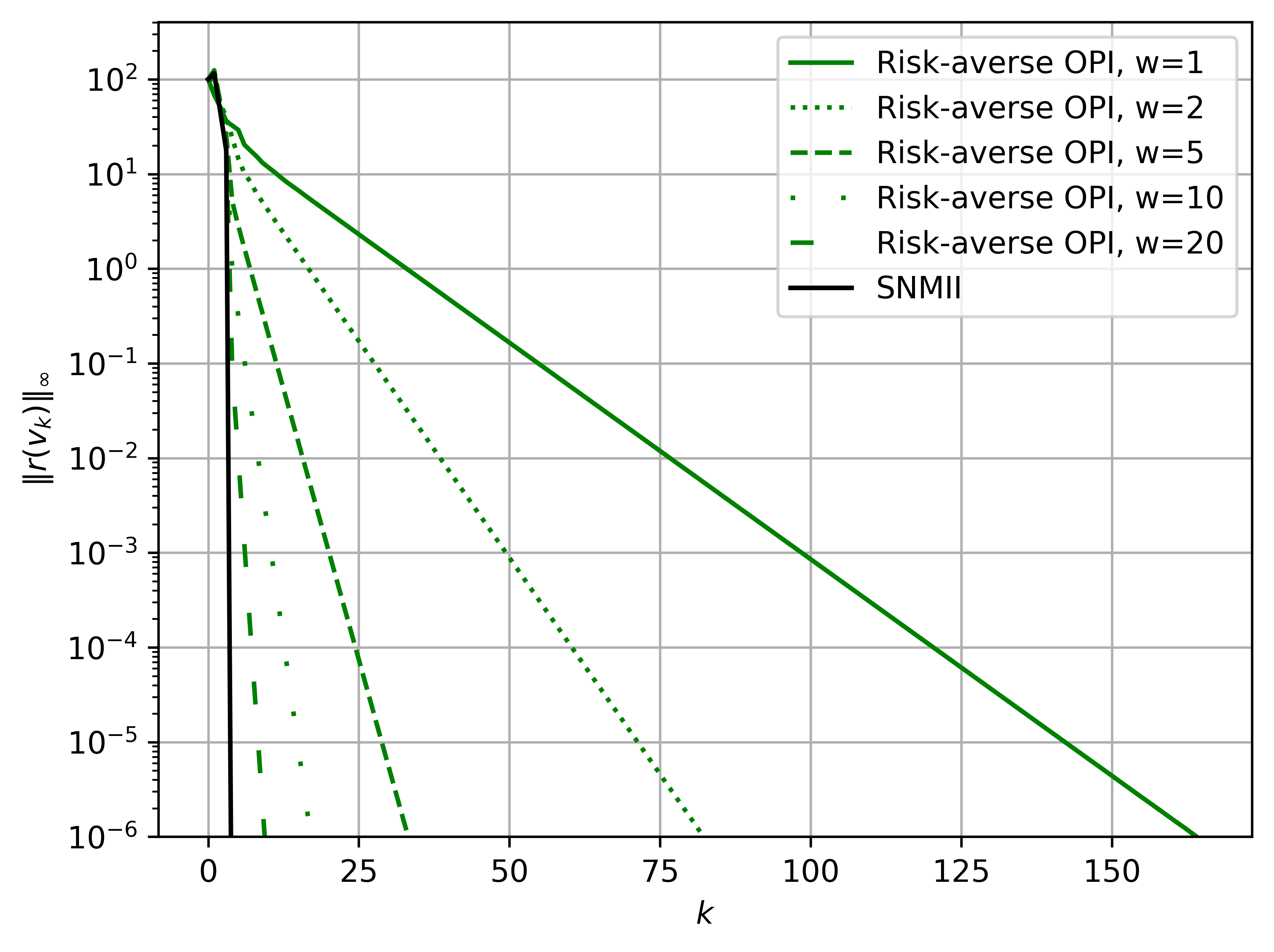}
    \caption{Infinity norm of the Bellman residual function versus the number of iterations for risk-averse OPI with a different number of inner iterations and SNMII.}
    \label{fig:RA_OPI}
\end{subfigure}
\caption{We consider an artificial MDP with $n=100$, $m=5$, $\gamma=0.9$, CVaR with $\zeta=0.3$ as risk-measure.}
\label{fig:riskaverseOPI}
\end{figure}

For the benchmarks we considered artificial MDPs where the transitions and stage-costs are generated by drawing samples from a uniform distribution. We tested the empirical performance of the methods on MDPs with different sizes in terms of state and action spaces, and we also considered a high ($\gamma=0.9$) and a low ($\gamma = 0.1$) discount factor. From our benchmarks we are able to conclude that the proposed SN methods have comparable performance in terms of both convergence rate (see results in Figure~\ref{fig:SNM}) and CPU time (see Table~\ref{table1_SNM} where we use SNMII$^*$ to denote SNMII with warm-start for the inner loop iterates).
In our experiments we noticed that SNMII and its warm-started variant tend to suffer from numerical instability. In particular, iteration~\eqref{eq:newton_inneriteration} fails in achieving the convergence when the inner tolerance level is set to a low value, \textit{e.g.}, $10^{-7}$, leading to oscillations. On the other hand, the convergence of SNMI and SNMIII is not compromised by the selection of the tolerance level parameters, which makes these methods more robust against the choice of their parameters setting.

Regarding risk-averse OPI, as to be expected from the theory, the higher the number of inner iterations executed, the faster the convergence rate in terms of number of outer iterations (see Figure~\ref{fig:RA_OPI}). Performing a higher number of inner iterations, however, comes with higher computational costs. As for risk-neutral OPI, the optimal number of inner iterations represents a trade-off between convergence rate and computational costs. In our experiments risk-averse OPI with a constant $w_k=20$ happens to be the best among the settings we studied. From Tables~\ref{table1_SNM} and~\ref{table2_OPI}, we conclude, however, that risk-averse OPI is not competitive 
with the other methods in terms of the overall computation time, unlike its risk neutral counterpart, even for a low discount factor. Indeed, value iteration for policy evaluation has a slow convergence rate in both the risk-neutral and risk-averse settings, but, while in the risk-neutral case its iterations are computationally cheap which makes the method competitive, in the risk-averse case the extra computational costs disrupt this trade-off. Finally, recall that risk-averse OPI with a constant $w_k=1$ is equivalent to value iteration. The results in Figure~\ref{fig:riskaverseOPI} clearly show the superior convergence rate of SN methods with respect to risk-averse value iteration. This is due to the faster local contraction properties of semismooth Newton methods with respect to risk-averse value iteration. In particular, as it is displayed in Figure~\ref{fig:SNM}, while SNMII (black bold line) reaches convergence to tolerance of $10^{-6}$ in less than 10 iterations, risk-averse value iteration (= risk-averse OPI, $w=1$ depicted with green bold line) requires more than 150 iterations to reach convergence to the same tolerance level. With respect to the risk-neutral scenario where slow rate of convergence is compensated by cheap iterations, the results of Table~\ref{table2_OPI} demonstrate that the extra computational costs arising from the evaluation of the risk-averse Bellman operator play against the risk-averse version of value iteration, which also in the risk-averse setting requires a substantially higher number of iterations than SN methods. These factors make the method not attractive. 

\begin{table}
\centering
\begin{tabular}{ |p{1.9cm}||p{1.1cm}|p{1.1cm}|p{1.1cm}|p{1.1cm}|  }
 \hline
 \multicolumn{5}{|c|}{\textbf{CPU Time [s]}} \\
 \hline
  $(n,\,m,\,\gamma)$ & SNMI & SNMII & SNMII$^*$ & SNMIII \\
 \hline
 $(50,\,5,\,0.9)$   & \textbf{2.71}    & 6.02 & 4.00 & 4.17 \\
  $(50,\,30,\,0.9)$ &   14.48  & 12.99   & \textbf{11.57} & 16.72 \\
 $(50,\,5,\,0.1)$ & 1.92 & 2.06 & \textbf{1.57} & 1.89 \\
 $(50,\,30,\,0.1)$  & 8.11 & 6.62 & \textbf{6.25} & 7.99 \\
 $(100,\,5,\,0.9)$ & \textbf{12.57} & 19.92 & 16.06 & 14.38\\
 $(100,\,20,\,0.9)$ & \textbf{41.91}  & 56.50 & 50.61 & 50.03\\
  $(100,\,5,\,0.1)$ & 7.35 & 8.52  & \textbf{6.88} & 7.30 \\
   $(100,\,20,\,0.1)$ & 27.99  & 24.14 & \textbf{22.73} & 28.94 \\
 \hline
\end{tabular}
\caption{We consider artificial MDPs with different size of state and action spaces, two different values of discount factor and CVaR with $\zeta=0.3$ as risk-measure. We report the CPU time in seconds for SNM methods and we use SNMII$^*$ to denote the warm-started version of SNMII.}
\label{table1_SNM}
\end{table}

\begin{table}
\centering
\begin{tabular}{ |p{1.9cm}||p{1.1cm}|p{1.1cm}|p{1.1cm}|p{1.1cm}|  }
 \hline
 \multicolumn{5}{|c|}{\textbf{CPU Time [s]}} \\
 \hline
  $(n,\,m,\,\gamma)$ &  OPI(2) & \!OPI(20) & \!OPI(50) & \!\!\!OPI(100) \\
 \hline
 $(50,\,5,\,0.9)$   & 37.02    & \textbf{15.38} & 19.19 & 36.80 \\
  $(50,\,30,\,0.9)$ &   194.22  & \textbf{36.26}   & 29.35 & 37.77 \\
 $(50,\,5,\,0.1)$ & \textbf{2.44} & 3.76 &  8.07 & 15.29 \\
 $(50,\,30,\,0.1)$  & 10.06 & \textbf{8.36} & 12.61 & 19.63 \\
 $(100,\,5,\,0.9)$ & 187.45 & \textbf{95.47} & 129.39 & 208.89 \\
 $(100,\,20,\,0.9)$ & 695.95 & \textbf{160.03} & 201.84 & 292.99 \\
  $(100,\,5,\,0.1)$ & \textbf{11.19}  & {21.65} & 46.13 & 87.96\\
   $(100,\,20,\,0.1)$ & 37.28 & \textbf{37.05} & 61.39 & 101.94\\
 \hline
\end{tabular}
\caption{For the same set of artificial MDPs of Table~\ref{table1_SNM}, we report the CPU time in seconds for risk-averse OPI with a different number of inner iterations. We denote with $\text{OPI}(i)$ the risk-averse OPI method with $i$-inner iterations.}
\label{table2_OPI}
\end{table}

\section{Conclusions \& Future Work}
We propose a new class of methods for risk-averse MDPs with finite state and action spaces. In particular, by exploiting the logic of semismooth Newton methods, we are able to derive three different methods to compute the root of the risk-averse Bellman residual function. One of the derived methods coincides with the risk-averse policy iteration method proposed by Ruszczyński in~\cite{risk_averse_DP}. Our semismooth Newton perspective not only allows one to formally explain the superior convergence properties of Ruszczyński's risk-averse policy iteration with respect to risk-averse value iteration, but it also offers a new framework for designing efficient methods to solve risk-averse MDPs.
As future work we plan to extend \textit{madupite}~\cite{madupite}, a high-performance distributed solver for risk-neutral MDPs, to also address the risk-averse setting by incorporating the proposed methods and relying on parallel/distributed computing to efficiently solve the convex programs.

\begin{ack}                               
This work was supported by the European Research Council under the Horizon 2020 Advanced under Grant 787845 (OCAL).
\end{ack}

\bibliographystyle{plain}        
\bibliography{autosam}           


\end{document}